\documentclass[12pt,leqno]{article} 
\usepackage{amssymb, amsmath, amsthm, amsxtra, subfigure, graphicx}
\newcommand{\ds}{\displaystyle}
\newcommand{\lf}{\left}
\newcommand{\rt}{\right}

\newcommand{\tit}{\textit}
\newcommand{\tbf}{\textbf}

\newcommand{\tdiag}{\text{diag}}
\newcommand{\tspec}{\text{spec}}

\newcommand{\beq}{\begin{equation}}
\newcommand{\eeq}{\end{equation}}

\newcommand{\llg}{\langle}
\newcommand{\rrg}{\rangle}

\newcommand{\RR}{\mathbb{R}}
\newcommand{\EE}{\mathbb{E}}
\newcommand{\mA}{\mathcal{A}}
\newcommand{\mE}{\mathcal{E}}
\newcommand{\tT}{\tilde{T}}
\newcommand{\ttt}{\tilde{t}}
\newcommand{\tX}{\tilde{X}}
\newcommand{\hk}{\hat{k}}

\newcommand{\De}{\Delta}
\newcommand{\ep}{\epsilon}

\newcommand{\Om}{\Omega}
\newcommand{\sg}{\sigma}
\newcommand{\Sg}{\Sigma}
\newcommand{\la}{\lambda}

\newcommand{\p}{\partial}

\newcommand{\half}{\frac{1}{2}}

\theoremstyle{plain}
\newtheorem{theorem}{Theorem}

\theoremstyle{definition}

\DeclareMathOperator{\prob}{Prob}

\begin{document}

\title{Universality in numerical computation with random data.  Case studies, analytic results and some speculations.}
\author{Percy Deift\thanks{The work in this paper was supported in part by 
DMS Grant \#13000965 (P.D.) and DMS Grant \#1303018 (T.T.).} \\
Courant Institute\\
\and
Thomas Trogdon$^*$\\
University of California, Irvine}
\date{}
\maketitle
\begin{abstract}
We discuss various universality aspects of numerical computations  using standard algorithms.  These aspects include empirical observations and rigorous results. We also make various speculations about computation in a broader sense.
\end{abstract}

\paragraph{Acknowledgements.} One of the authors (P.D.) would like to thank the organizers for the invitation to participate in  the Abel Symposium 2016 ``Computation and Combinatorics in Dynamics, Stochastics and Control". During the symposium he gave a talk on a condensed version of the paper below.

\newpage

There are two natural ``integrabilities'' associated with matrices $M$.  The
first concerns random matrix theory where key statistics, such as the 
distribution of the largest eigenvalue of $M$, or the gap probability, i.e.,
the probability  that the spectrum of $M$ contains a prescribed gap, are
described in an appropriate scaling limit as $N=\dim M \to \infty$, by the
solution of completely integrable Hamiltonian systems, viz., Painlev\'e
equations (see e.g.\ \cite{Meh}).  The second concerns the numerical computation of the eigenvalues of a matrix. Standard eigenvalue algorithms work in the 
following way.  Let $\Sg_N$ denote the set of real $N\times N$ symmetric
matrices and let $M\in \Sg_N$ be a given matrix whose eigenvalues one wants
to compute. Associated with each algorithm $\mA$, there is, in the discrete
case, a map $\varphi = \varphi_\mA:\Sg_N\to \Sg_N$, with the properties
\begin{itemize}
\item (isospectral)\qquad  $\tspec \lf(\varphi_\mA (H)\rt) = 
\tspec \; (H), \qquad H \in \Sg_N$,
\item (convergence) \quad the iterates $X_{k+1} = \varphi_\mA (X_k),
\quad k\ge 0, \qquad X_0=M$, converge to a diagonal matrix $X_\infty, \quad X_k \to X_\infty$,\quad as
\quad $k\to\infty$,
\end{itemize}
and in the continuum case, there is a flow $t\to X(t)
\in \Sg_N$ with the properties
\begin{itemize}
\item (isospectral) $\quad\tspec \lf(X(t)\rt) = \tspec\, \lf(X(0)\rt)$,
\item (convergence) the flow $X(t)$, $t\ge 0$, $X(0)=M$, converges to a 
diagonal matrix $X_\infty$, $X(t)\to X_\infty$ as $t\to \infty$.
\end{itemize}
In both case, necessarily the (diagonal) entries of $X_\infty$ are the 
eigenvalues  of the given matrix $M$.  Now the fact of the matter is that,
in most cases of interest, the flow $t\to X(t)$ is Hamiltonian and completely
integrable in the sense of Liouville, and in the discrete case we have a
``stroboscope theorem'', i.e.\ there exists a completely integrable 
Hamiltonian flow $t\to \tX(t)$ which coincides with the above iterates $X_k$
at integer times, $\tX(k)=X_k, \;\; k\ge 0$ (see, in particular, 
\cite{Sym}, \cite{DNT}, \cite{DLNT}).  The abstract QR algorithm 
is a prime example of such a discrete
algorithm, while the Toda algorithm is an example of the continuous case.

Question: What happens if one tries to ``marry'' these two integrabilities?
In particular, what happens when one computes the eigenvalues of a random
matrix?  In response to this question, the authors in \cite{PDM} initiated a 
statistical study of the performance of various standard algorithms to compute
the eigenvalues of random matrices $M$ from $\Sg_N$.

Given $\ep >0$, it follows, in the discrete case, that for some $m$ the 
off-diagonal entries of $X_m$ are $O(\ep)$ and hence the diagonal entries 
of $X_m$ give the eigenvalues of $X_0=M$ to $O(\ep)$.  The situation is 
similar for continuous algorithms $t\to X(t)$.  Rather than running the
algorithm until all the off-diagonal entries are $O(\ep)$, it is customary 
to run the algorithm with \tbf{deflations} as follows.  For an $N\times N$
matrix $Y$ in block from 
$$
Y= \lf[\begin{matrix}
Y_{11} & Y_{12}\\
Y_{21} &Y_{22}
\end{matrix}
\rt]
$$
with $Y_{11}$ of size $k\times k$ and $Y_{22}$ of size $(N-k) \times (N-k)$
for some $k\in \lf\{1, 2, \dots, N-1\rt\}$, the process of projecting $Y\to
\tdiag \; \lf(Y_{11}, Y_{22}\rt)$ is called deflation.   For a given
$\ep >0$, algorithm $\mA$ and matrix $M\in \Sg_N$, define the 
{\boldmath $k$}\tbf{-deflation} time $T^{(k)}(M) = T^{(k)}_{\ep, \mA}(M),
\quad 1\le k \le N-1$, to be the smallest value of $m$ such that $X_m$,
the $m^{th}$ iterate of algorithm $\mA$ with $X_0 = M$, has block form
$$
X_m= \lf[\begin{matrix}
X^{(k)}_{11} & X^{(k)}_{12}\\
X^{(k)}_{21} &X^{(k)}_{22}
\end{matrix}
\rt]
$$
with $X^{(k)}_{11}$ of size $k\times k$ and $X^{(k)}_{22}$ of size 
$(N-k) \times (N-k)$ and $\| X^{(k)}_{12}\| = \|X^{(k)}_{21}\| \le \ep$.
The deflation time $T(M)$ is then defined as
$$
T(M) = T_{\ep,\mA}(M) = \min_{1\le k \le N-1} \quad T^{(k)}_{\ep,\mA} (M).
$$
If $\hk\in \lf\{1, \dots, N-1\rt\}$ is such that $T(M)= T^{(\hk)}_{\ep, \mA}
(M)$, it follows that the eigenvalues of $M=X_0$ are given by the eigenvalues
of the block-diagonal matrix $\tdiag \lf(X^{(\hk)}_{11}, X^{(\hk)}_{22}\rt)$
to $O(\ep)$.  After running the algorithm to time $T_{\ep, \mA}(M)$, the
algorithm restarts by applying the basic algorithm $\mA$ separately to the
smaller matrices $X^{(\hk)}_{11}$ and $X^{(\hk)}_{22}$ until the next deflation
time, and so on.  There are again similar considerations for continuous 
algorithms.  

As the algorithm proceeds, the number of matrices after each deflation doubles.
This is counterbalanced by the fact that the matrices are smaller and smaller 
in size, and the calculations are clearly parallelizable.  Allowing for parallel
computation, the number of deflations to compute all the eigenvalues of a given 
matrix $M$ to an accuracy $\ep$, will vary from $O(\log N)$ to $O(N)$.

\begin{figure}[htbp]
\centering
\subfigure[]{\includegraphics[width=.45\linewidth]{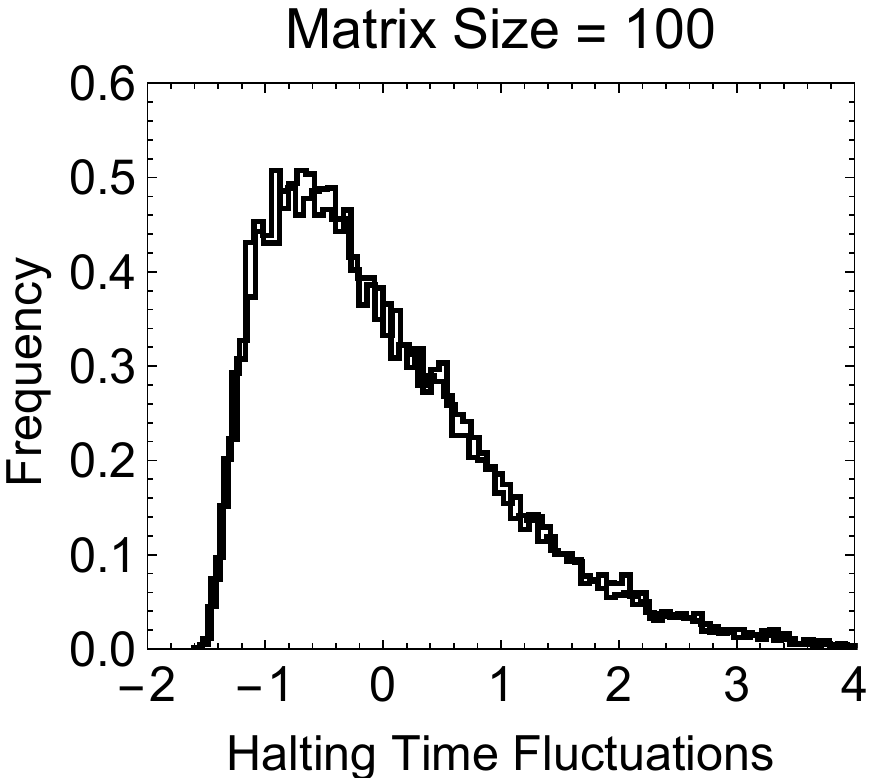}}
\subfigure[]{\includegraphics[width=.45\linewidth]{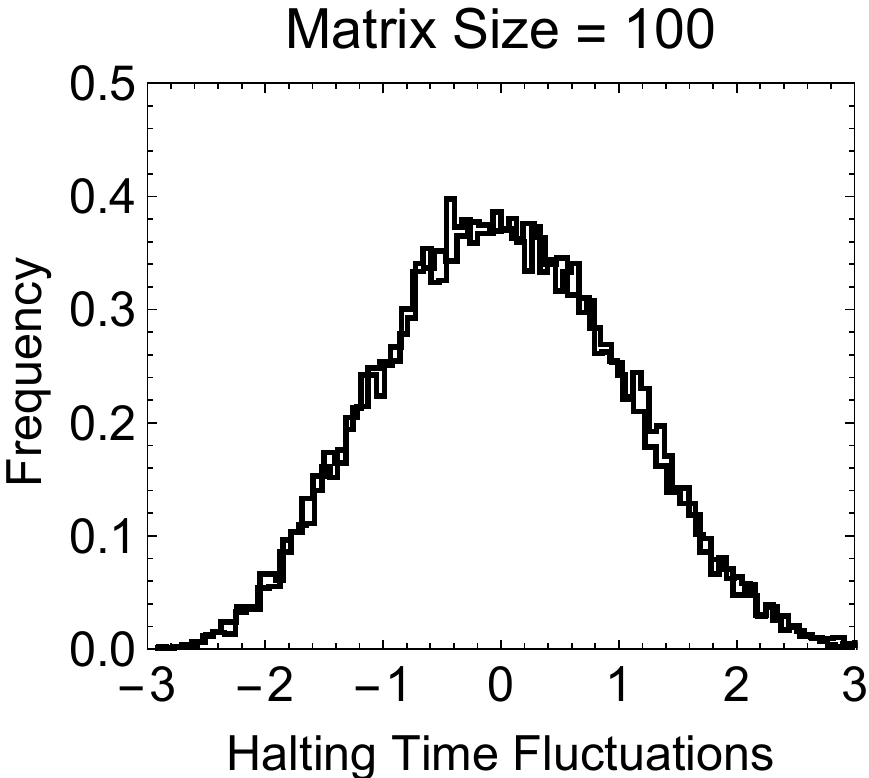}}
\caption{\label{f:one} Universality for $\tilde T$ when (a) $\mathcal A$ is the QR eigenvalue algorithm and when (b) $\mathcal A$ is the Toda algorithm.  Panel (a) displays the overlay of two histograms for $\tilde T$ in the case of QR, one for each of the two ensembles $ \mathcal E = \mathrm{BE}$, consisting of iid mean-zero Bernoulli random variables and $\mathcal E = \mathrm{GOE}$, consisting of iid mean-zero normal random variables. Here $\epsilon= 10^{-10}$ and $N = 100$.  Panel (b) displays the overlay of two histograms for $\tilde T$ in the case of the Toda algorithm, and again $\mathcal E = \mathrm{BE}$ or $\mathrm{GOE}$. And here $\epsilon= 10^{-8}$ and $N = 100$. }
\end{figure}

In \cite{PDM} the authors considered the deflation time $T=T_{\ep,\mA}=
T_{\ep, \mA, \mE}$ for $N\times N$ matrices chosen from an ensemble $\mE$.
For a given $\ep >0$, algorithm $\mA$ and ensemble $\mE$, the authors computed
$T(M)$ for 5,000--10,000 samples of matrices $M$ chosen from $\mE$, and
recorded the \tbf{normalized deflation time}
\beq\label{1}
\tT(M) \equiv \frac{T(M)-\llg T\rrg}{\sg}
\eeq
where $\llg T\rrg$ and $\sg^2= \lf\llg \lf(T-\llg T\rrg\rt)^2\rt\rrg$
are the sample average and sample variance of $T(M)$, respectively.  What the 
authors found, surprisingly, was that for the given algorithm $\mA$, and
$\ep$ and $N$ in a suitable scaling range with $N\to \infty$, the 
\tbf{histogram of} {\boldmath $\tT$} \tbf{was universal, independent of the
ensemble } {\boldmath $\mE$}.  In other words, the fluctuations in the 
deflation time $\tT$, suitably scaled, were universal, independent of $\mE$.
Figure~\ref{f:one} displays some of the numerical results from \cite{PDM}.
Figure 1(a) displays data for the QR algorithm, which is discrete, and 
Figure 1(b) displays data for the Toda algorithm, which is continuous.  
Note that the histograms in Figures 1(a) and 1(b) are very different:
Universality is observed with respect to the ensembles $\mE$---not with respect
to the algorithms $\mA$. 

Subsequently in \cite{DMOT} the authors raised the question of whether the
universality results in \cite{PDM} were limited to eigenvalue algorithms 
for real symmetric matrices, or whether they were present more generally 
in numerical computation.  And indeed the authors in 
\cite{DMOT} found similar 
universality results for a wide variety of numerical algorithms, including
\begin{enumerate}
\item[(a)] other algorithms such as the QR algorithm with shifts\footnote{The QR algorithm with shifts is the accelerated version of the QR algorithm that is used in practice.}, the Jacobi
eigenvalue algorithm, and also algorithms applied to complex Hermitian 
ensembles
\item[(b)] the conjugate gradient and GMRES algorithms to solve linear
$N\times N$ systems $Hx=b$ with $H$ and $b$ random
\item[(c)] an iterative algorithm to solve the Dirichlet problem $\De u=0$
in a random star-shaped region $\Om \subset \RR^2$ with random boundary data
$f$ on $\p \Om$
\item[(d)] a genetic algorithm to compute the equilibrium measure for
orthogonal polynomials on the line.
\end{enumerate}

In \cite{DMOT} the authors also discuss similar universality results 
obtained by Bakhtin and Correll \cite{BC} in a series of experiments with 
live participants
recording 
\begin{enumerate}
\item[(e)] decision making times for a specified task.
\end{enumerate}
Whereas (a) and
(b) concern finite dimensional problems, (c) shows that universality is also
present in problems that are genuinely infinite dimensional.  And whereas
(a), (b) and (c) concern, in effect, deterministic dynamical systems acting
on random initial data, problem (d) shows that universality is also present
in genuinely stochastic algorithms.

The demonstration of universality in problems (a)--(d) raises the following
issue: Given the common view of neuroscientists that the brain is just a big
computer with hardware and software, one should be able to find evidence of
universality in some neural computations.  It is this issue that led the authors
in \cite{DMOT} to the work of Bakhtin and Correll.  In \cite{BC} 
each of the participants
is shown a large number $k$ of  diagrams and then asked to make a decision 
about a particular geometric feature of each diagram.  What is then recorded
is the time it takes for the participant to reach his'r decision.  Thus each
participant produces $k$ decision times $t$ which are then centered and scaled
as in \eqref{1} to obtain a normalized decision time
\beq\label{2}
\ttt = \frac{t-\llg t \rrg}{\sg}.
\eeq
The distribution of $\ttt$ is then recorded in a histogram.  Each of the 
participants produces such a histogram, and what is remarkable is that the 
histograms are, with a few exceptions, (essentially) the same.  Furthermore,
in \cite{BC}, Bakhtin and Correll developed a Curie-Weiss-type statistical
mechanical model for the decision process, and obtained a distribution function
$f_{BC}$ which agrees remarkably well with the (common) histogram obtained 
by the participants.  We note that the model of Bakhtin and Correll involves
a particular parameter, the spin flip intensity $c_i$.  In \cite{BC} 
the authors made one particular choice for $c_i$.  However, as shown in 
\cite{DMOT}, if one
makes various other choices for $c_i$, then one still obtains the same
distribution $f_{BC}$.  In other words, the Bakhtin-Correll model itself has
an intrinsic universality.  In an independent development Sagun, Trogdon and 
LeCun \cite{STL} considered, amongst other things, search times on {\tt Google$^\mathrm{TM}$} for a 
large number of words in English and in Turkish.  They then centered and scaled
these times as in \eqref{1}, \eqref{2} to obtain two histograms for normalized search times,
one for English words and one for Turkish words.  To their great surprise, both
histograms were the same and, moreover, extremely well described by $f_{BC}$.
So we are left to ponder the following puzzlement: Whatever the neural
stochastics of the participants in the study in \cite{BC}, and whatever the 
stochastics in the Curie-Weiss model, and whatever the mechanism in {\tt Google$^\mathrm{TM}$}'s
search engine, a commonality is present in all three cases expressed through
the single distribution function $f_{BC}$.

All of the above results are numerical.  In order to establish universality 
as a bona fide phenomenon in numerical analysis, and not just an artifact,
suggested, however strongly, by certain computations as above, P.~Deift
and T.~Trogdon in \cite{DT1} considered the Toda eigenvalue algorithm mentioned
above.  In place of the deflation time $T(M)=\min_{1\le k \le N-1}\;\;
T^{(k)}_{\ep\,\mA}(M)$, $\mA=\text{Toda algorithm}$, Deift and Trogdon used
the 1-deflation time $T^{(1)}(M)= T^{(1)}_{\ep, \mA}(M)$ as the stopping time
for the algorithm.  In other words, given $\ep >0$ and an ensemble $\mE$,
they ran the Toda algorithm $t\to X(t)$ with $X(0)=M \in \mE$, until a time
$t$ where
$$
t=T^{(1)}(M) = \inf \lf\{s\ge 0 : \sum^N_{j=2} \lf(X_{1j}
(s)\rt)^2 \le \ep^2\rt\}.
$$
It follows by perturbation theory that $\lf|X_{11} \lf(T^{(1)}(M)\rt)
- \la_{j^\ast}(M)\rt| \le \ep$ for some eigenvalue $\la_{j^\ast} (M)$
of $M$. But the Toda algorithm is known to be ordering, i.e.\
$X(t)\to X_\infty=\tdiag \lf(\la_1(M), \la_2(M), \dots \la_N(M)\rt)$, 
where the eigenvalues of $M$ are ordered, $\la_1(M) \ge \la_2(M) \ge
\dots \ge \la_N(M)$.  It follows then that (for $\ep$ sufficiently small
and $T^{(0)}_{\ep, \mA}$ correspondingly large) $j^{\ast} =1$ so that
the Toda algorithm with stopping time $T^{(1)}= T^{(1)}_{\ep, \mA}$
computes the largest eigenvalue of $M$ to accuracy $\ep$ with high probability.

The main result in \cite{DT1} is the following.  For invariant and generalized 
Wigner random matrix ensembles\footnote{See Appendix A in \cite{DT1} 
for a precise
description of the matrix ensembles considered in Theorem 1.} there is 
an ensemble dependent constant $c_{\mE}$ such that the  following limit
exists (see \cite{PS} and \cite{WBF})
\beq\label{3}
F^{\text{gap}}_{\beta} (t) = \lim_{N\to\infty} \prob \lf(
\frac{1}{c^{2/3}_{\mE} \; 2^{-2/3}\; N^{2/3} \lf(\la_1-\la_{2}\rt)} 
\le t\rt), \qquad t\ge 0.
\eeq
Here $\beta=1$ for the real symmetric case, 
$\beta=2$ for the complex Hermitian case. Thus $F^{\text{gap}}_\beta(t)$ is the distribution function for the 
(inverse of the) gap $\la_1-\la_{2}$ between the largest two eigenvalues
of $M$, on the appropriate scale as $N\to \infty$.

\begin{theorem}[Universality for $T^{(1)}$]
Let $0< \sg <1$ be fixed and let $(\ep, N)$ be in the scaling region 
\beq\label{4}
\frac{\log \ep^{-1}}{\log N} \ge \frac{5}{3} + \frac{\sg}{2}\ .
\eeq
Then if $M$ is distributed according to any real $(\beta=1$) or complex
$(\beta=2)$ invariant or Wigner ensemble, we have 
\beq\label{5}
\lim_{N\to \infty} \prob \lf( 
\frac{T^{(1)}}{ c^{2/3}_{\mE}\;2^{-2/3}\;N^{2/3} \lf(\log \ep^{-1}-
\frac{2}{3} \;\log N\rt) } \le t \rt)= F^{\text{gap}}_\beta.
\eeq
Here $c_\mE$ is the same constant as in (3).
\end{theorem}

This result establishes universality rigorously for a numerical algorithm
of interest, viz., the Toda algorithm with stopping time $T^{(1)}$ to compute
the largest eigenvalue of a random matrix.  We see, in particular, that 
$T^{(1)}$ behaves statistically as the inverse of the top gap 
$\la_1 - \la_{2}$, on the appropriate scale as $N\to \infty$.  Similar
results have now been obtained for the QR algorithm and related algorithms
acting on ensembles of strictly positive definite matrices (see \cite{DT2}).

The proof of Theorem 1 depends critically on the integrability of the Toda
flow $t\to X(t),\;\;X(0)=M$. The evolution of $X(t)$ is governed by the
Lax-pair equation 
$$
\frac{dX}{dt} = \lf[ X, \, B(X)\rt] = X\, B(X)- B(X)\,X
$$
where $B(X) = X_- - X^T_-$ and $X_-$ is the strictly lower triangular part of $X$.
Using results of J. Moser \cite{Mos} one finds that
\begin{align}
&E(t) \equiv \sum^N_{k=2} \lf|X_{1\,k}(t)\rt|^2 = \sum^N_{j=1}
\lf(\la_j-X_{11}(t)\rt)^2 \lf|u_{1j}(t)\rt|^2 \label{6}\\
&\hspace{.10in}X_{11}(t) = \sum^N_{j=1} \la_j \lf|u_{1j}(t)\rt|^2 \label{7} \\
&u_{1j}(t) = \frac{ u_{1j}(0)\, e^{\la_j\, t} }{\lf( \ds{\sum^N_{k=1}} 
\lf| u_{1k}(0)\rt|^2\, e^{2\la_k\, t}\rt)^{\half} }, \qquad 1\le j\le N,
\label{8}
\end{align}
where $u_{1j}(t)$ is the first component of the normalized eigenvector $u_j(t)$
for $X(t)$ corresponding  to the eigenvalue $\la_j(t) = \lambda_j(0)$ of $X(t)$, 
$\lf(X(t)-\la_j(t)\rt) u_j(t)=0$.  (Note that $t\to X(t)$ is isospectral, 
so $\tspec(X(t)) = \tspec(X(0)) = \tspec(M)$.)  The stopping time $T^{(1)}$ is obtained
by solving the equation 
\beq\label{9}
E(t) = \ep^2
\eeq
for $t$.  Substituting \eqref{7} and \eqref{8} into \eqref{6} we obtain an formula for $E(t)$
involving only the eigenvalues and (the moduli of) the first components of 
the normalized eigenvectors for $X(0)=M$.  It is this explicit formula that
the Toda algorithm brings as a gift to the marriage announced earlier of
eigenvalue algorithms and random matrices.  What random matrix theory brings
to the marriage  is an impressive collection of very detailed estimates on
the statistics of the $\la_j$'s and the $u_{1j}(0)$'s obtained in recent
years by a veritable army of researchers including P.~Bourgade, L.~Erd\H{o}s,
A.~Knowles, J. A. Ram\'irez, B. Rider, B. Vi\'rag, T. Tao, V. Vu, J. Yin
and H. T. Yau, amongst many others (see \cite{DT1} and the references therein 
for more details).

Theorem 1 is a first step towards proving universality  for the Toda algorithm
with full deflation stopping time $T=T_{\ep, \mA}$.  The analysis of $T_{\ep,
\mA}$ involves very detailed information about the joint statistics of the 
eigenvalues $\la_j$ and all the components $u_{ij}$ of the normalized 
eigenvectors of $X(0)=M$, as $N\to\infty$.  Such information is not yet 
known and the analysis of $T_{\ep, \mA}$ is currently out of reach.

\paragraph{Speculations.} How should one view the various \tbf{two-component} universality results 
described in this paper?  ``Two-components'' refers to the fact for a random
system of size $S$, say, and halting time $T$, once the average $\llg T\rrg$
and variance $\sg^2= \lf\llg \lf(T - \llg T\rrg \rt)^2\rt\rrg$ are known,
the normalized time $\tau=\lf(T-\llg T\rrg\rt)/\sg$ is, in the large $S$
limit, universal, independent of the ensemble, i.e. as $S\to\infty$, $T\sim
\llg T\rrg + \sg\, \chi$, where $\chi$ is universal.  The best known two-component
universality theorem is certainly the classical Central Limit Theorem:
Suppose $Y_1, Y_2, \dots$ are independent, identically distributed variables
with mean $\mu$ and variance $\sg^2$. Set $W_n\equiv \sum^n_{i=1}\, Y_i$.
Then as $n\to\infty$, $\lf(W_n-\llg W_n\rrg\rt)/\sg_n$ converges in 
distribution to a standard normal $N(0, 1)$, where $\llg W_n\rrg = \EE \left( \sum_{i=1}^n Y_i \right) = n \mu$ and $ \sigma_n^2 = \EE
\lf(\lf(W_n- \llg W_n\rrg\rt)^2\rt)= n\sg^2$.  In words: As $n\to\infty$,
the only specific information about the initial distribution of the $Y_i$'s
that remains, is the mean $\mu$ and the variance $\sg$.

Now imagine you are walking on the boardwalk in some seaside town.  Along 
the way you pass many palm trees.  But what do you mean by a ``palm tree''?
Some are taller, some are shorter, some are bushier, some are less bushy.
Nevertheless you recognize them all as ``palm trees'': Somehow you adjust
for the height and you adjust for the bushiness (two components!), and 
then draw on some internal data base to determine, with high certainty,
that the object one is looking at is a ``palm tree''.  The database itself catalogs/summarizes your learning experience with palm trees over many years. It is tempting to
speculate that the data base has the form of a histogram.  We have in our
brains one histogram for palm trees, and  another for olive trees, and so on.
Then just as we may use a $t$-test, for example, to test the statistical 
properties of some sample, so too one speculates that there is a mechanism in 
one's mind that tests against the ``palm tree histogram'' and evaluates the 
likelihood that the object at hand is a palm tree.  So in this way of thinking,
there is no ideal Platonic object that is a ``palm tree'':  Rather, a palm
tree is a histogram.

One may speculate further in the following way.  Just imagine if we perceived
every palm tree as a distinct species, and then every olive tree as a distinct
species, and so on. Working with such a plethora of data, would require access
to an enormous bandwidth.  From this point of view, the histogram provides a
form of ``stochastic data reduction'', and the fortunate fact is that we have 
evolved to the point that we have just enough bandwidth to accommodate and 
evaluate the information ``zipped'' into the histogram.  On the other hand,
fortunately, the information in the histogram is sufficiently detailed that we
can make meaningful distinctions, and one may speculate that it is precisely
this balance between data reduction and bandwidth that is the key to our
ability to function successfully in the macroscopic world.

We note finally that there are many similarities between the above speculations and machine learning.  In both processes there is a learning phase followed by a recognition phase.  Also, in both cases, there is a balance between data reduction and bandwidth. In the case of the palm trees, etc., however we make the additional assertion/speculation that the stored data is in the form of a histogram, similar in origin to the universal histograms observed in numerical computations.

We may summarize the above discussion and speculations in the following way.
The brain is a computer, with software and hardware, which makes calculations and 
runs algorithms which reduce data on an appropriate scale---the
macroscopic scale on which we live---to a manageable and useful form, viz.,
a histogram, which is universal\footnote{A priori the histogram for a palm tree in one person's mind may be very different from that in another person's mind.  Yet the results of Bakhtin and Correll in \cite{BC}, where the participants produce the same decision time distributions, indicate that this is not so.  And indeed, if there was a way to show that the histograms individuals form to catalog a palm tree, say, were all the same, this would have the following implication: The palm tree has an objective existence, and not a subjective one, which varies from person to person.} for all palm trees, or all olive trees, etc. 
With this in mind, it is tempting to suggest that \tbf{whenever we run an 
algorithm with random data on a ``computer", two-component universal features will emerge
on some appropriate scale.} This ``computer'' could be the electronic machine on our desk, or it could be the device in our mind that runs algorithms to classify random visual objects or to make timed decisions about geometric shapes, or it could be in any of the myriad of ways in which computations are made. Perhaps this is how one should view the various
universality results described in this paper.

\end{document}